\input amstex 
\documentstyle{amsppt}
\input bull-ppt
\keyedby{bull428e/PAZ}


\define\br{{\Bbb R}}
\define\bz{{\Bbb Z}}
\define\cB{{\Cal B}}
\define\cC{{\Cal C}}
\define\cO{{\Cal O}}
\define\cL{{\Cal L}}

\define\sqr#1#2{{\vcenter{\vbox{\hrule height.#2pt
  \hbox{\vrule width.#2pt height#1pt \kern#1pt
  \vrule width.#2pt}
  \hrule height.#2pt}}}}


\hyphenation{meas-ured}
\hyphenation{meas-ure-ment}
\hyphenation{tatonne-ment}

\topmatter
\cvol{29}
\cvolyear{1993}
\cmonth{October}
\cyear{1993}
\cvolno{2}
\cpgs{228-234}
\title Harmonic Analysis of Fractal Measures Induced\\ 
by Representations of a Certain C$^*$-Algebra \endtitle
\shorttitle{Harmonic analysis of fractal measures}
\author Palle E. T. Jorgensen and Steen Pedersen \endauthor
\shortauthor{P. E. T. Jorgensen and Steen Pedersen}
\address Department of Mathematics, The University of 
Iowa, Iowa City, Iowa 
52242\endaddress
\ml jorgen\@math.uiowa.edu\endml
\address Department of Mathematics, Wright State 
University, Dayton, Ohio  45435\endaddress
\ml spedersen\@desire.wright.edu\endml
\thanks Research supported by the NSF.
The first author was partially supported
by a University of Iowa Faculty Scholar Award and the UI
(Oakdale Campus) Institute for Advanced Studies\endthanks
\date September 25, 1992\enddate
\abstract We describe a class of measurable subsets 
$\Omega$ in $\br^d$ such
that $L^2(\Omega)$ has an orthogonal basis of frequencies
$e_\lambda(x)=e^{i2\pi\lambda\cdot x}(x\in\Omega)$ indexed 
by
$\lambda\in\Lambda\subset\br^d$.  We show that such 
spectral pairs $(\Omega
,\Lambda)$ have a self-similarity which may be used to 
generate associated
fractal measures $\mu$ with Cantor set support.  The 
Hilbert space $L^2(\mu)$
does not have a total set of orthogonal frequencies, but a 
harmonic analysis of
$\mu$ may be built instead from a natural representation 
of the Cuntz C$^*$-
algebra which is constructed from a pair of lattices 
supporting the given
spectral pair $(\Omega ,\Lambda)$.  We show conversely 
that such a pair may be
reconstructed from a certain Cuntz-representation given to 
act on $L^2(\mu)$.\endabstract
\endtopmatter

\document

\heading 1.  Introduction\endheading

Let $\Omega$ be a subset in $d$ real dimensions (i.e., 
$\Omega\subset\br^d$,
$d\geq 1$), and suppose that $\Omega$ has finite positive 
$d$-dimensional
Lebesgue measure.  Let $L^2(\Omega)$ be the corresponding 
Hilbert space with
the usual inner product given by $$ \langle f,g\rangle =
m_d(\Omega)^{-1}\int_\Omega \overline{f(x)}\ g(x)\ dx$$ 
where $dx:=dx_1\cdots
dx_d$, and $m_d(\Omega)$ denoting the Lebesgue measure of 
$\Omega$.  Motivated
by a problem of I. E. Segal  and a paper by B. Fuglede 
[Fu], we considered in
[JP1--3] the problem of deciding, for given $\Omega$, when 
$L^2(\Omega)$ may
possibly have an orthogonal basis of frequencies:  For 
$\lambda\in\br^d$, let
$x\cdot\lambda =\sum^d_{j=1} x_j\lambda_j$ be the usual 
dot product, and set
$$e_\lambda(x)=e^{i2\pi\ x\cdot\lambda}\ .\tag 1$$ We say 
that two vector
frequencies $\lambda ,\lambda'$ in $\br^d$ are {\it 
orthogonal} on $\Omega$
if
$$\int_\Omega e^{i2\pi (\lambda'-\lambda)\cdot x} dx=0\ .$$
When $\Omega$ is further assumed open in $\br^d$, this 
problem is directly
connected (see  [Fu, JP1]) with the problem of finding 
simultaneous commuting
selfadjoint extension operators for the partial 
derivatives $\sqrt{-1}
{\partial\over\partial x_j}$   $(1\leq j\leq d)$ acting on 
$C^\infty_c
(\Omega)$ (= all smooth compactly supported functions in 
$\Omega)$.  In
general, the problem may be given a group-theoretic 
\pagebreak formulation, and, in this
form, we showed in [JP1] that it relates directly to a 
property of the
representation ring generated by a certain induced 
representation. (See (2)
below.)

\heading 2.  Classical Examples\endheading

The most obvious examples of sets $\Omega$ with the basis 
property are measurable sets in $\br^d$ which are {\it 
fundamental domains} of lattices (see
[Fu, JP1]).  Let $\Gamma$ be a rank $d$ {\it lattice}, and 
let $\Gamma\,^0$ be the {\it dual lattice}.
$$
(\text{Recall}\ \ \ \Gamma\,^0=\{ \lambda\in\br^d:\lambda
\cdot s\in\bz\ ,\ \forall s\in\Gamma\}.)$$
Suppose $\Omega$ is a measurable fundamental domain for 
$\Gamma$.  It is a simple matter to show then
that $\{e_\lambda :\lambda\in\Gamma\,^0\}$ is an
orthogonal basis for $L^2(\Omega)$.  This elementary class 
of examples is in fact characterized by a multiplicative 
property 
(see [JP1, 2]), and they
are called {\it multiplicative}.  A pair---$(\Omega 
,\Lambda)$ such that $0\in\Lambda$, and $\{ e_\lambda 
:\lambda\in\Lambda\}$ is an orthogonal basis
in $L^2(\Omega)$---is called a {\it spectral pair}, and 
the set $\Lambda$ is called the {\it spectrum}.  We 
further showed in [JP1] that every spectral
pair $(\Omega ,\Lambda)$ in $d$ dimensions may be 
factored, $(\Omega
,\Lambda)\simeq (\Omega',\Lambda')\times 
(\Omega'',\Lambda'')$,
 such that the
factors each are spectral pairs in dimensions $d', d''$ 
respectively, $d'+d''=d$, $(\Omega',\Lambda')$ is 
multiplicative, and $(\Omega'',\Lambda'')$
is in ``the other extreme''.  Specifically, this second 
factor generates a representation ring which is a copy of 
the algebra of all $q$ by $q$
complex matrices where $q$ is a certain cover-multiplicity 
(see [JP2]), and $(\Omega'',\Lambda'')$ is called a {\it 
simple factor}.

\heading 3.  Spectral pairs\endheading

In this paper, we shall consider the simple factors in 
more detail and show that they are associated with 
``fractals'' in a sense which we proceed to
describe.  If $(\Omega ,\Lambda)$ is a spectral pair in 
$d$ dimensions, consider the group 
$K=\Lambda^0=\{s\in\br^d:s\cdot\lambda\in\bz$,
$\forall\lambda \in\Lambda\}$.  We further showed in [JP1] 
that $K$ is a rank $d$ lattice and that there is a 
canonical embedding of $\Omega$ into
the torus $\br^d/K$ such that the image $\Omega'$ of 
$\Omega$ on the torus again has the basis-property 
(relative to Haar measure on the 
torus) and
the spectrum of $\Omega'$ is the same set $\Lambda$.  We 
say that the pair $(\Omega',\Lambda)$ is in {\it reduced 
form}.

We have a second closed subgroup $A$ in $\br^d$ directly 
associated with some given spectral pair $(\Omega 
,\Lambda)$,
$$A=\{ a\in\br^d:x+a\in\Omega +\Lambda^0 \roman{\ (a.e.)\ 
}x\in\Omega\}\ .$$
Define a unitary representation $U_t\ (t\in\br^d)$, acting 
on $L^2(\Omega)$, given by
$$U_t e_\lambda =e^{i2\pi t\cdot\lambda} e_\lambda\qquad 
(t\in\br^d,\ \lambda\in\Lambda)\ ,\tag 2$$
and note that $A$ may be characterized alternatively as 
the group
$$\{ t\in\br^d:U_t\hskip10pt\text{acts multiplicatively 
on}\hskip10pt L^2(\Omega)\}\ .$$
When $t\in A$, then
$$U_tf(x)=f(x+t)\ ,\quad \text{a.e.}\ \ x\in\Omega'\ ,\ 
\forall f\in L^2(\Omega')\tag 3$$
where the sum $x+t$ is in the torus $\br^d/K$.  Hence, we 
get $A$ acting as a group of torus-translations on 
$\Omega'$.

We say that some given spectral pair $(\Omega ,\Lambda)$ 
is {\it multiplicative} 
if $A=\br^d$ and is a {\it simple factor} if $A$ is a 
lattice in
$\br^d$.  There is a sense in which simple factors may be 
generated by lattice systems, but we do not yet have a 
complete structure theorem which
covers all simple factors.  It is not known if, for a 
simple factor with associated lattices $K$ and $A$, the 
{\it degenerate} case $K=A$ may occur.
(We expect not!)  In [JP1], we proved the following result 
(which will be needed below) about nondegenerate simple 
factors:

\topspace{25.5pc}\caption{Figure 1.}

\proclaim{Theorem 1 \rm(see [JP1], Theorem 6.1])}
 Let $(\Omega ,\Lambda)$ be a spectral pair in $\br^d$, 
and suppose that the group $S$, given by
$S=\{s\in\Lambda : s+\Lambda =\Lambda\}$, is a lattice.  
Let $\Gamma =S^0$, and suppose
\roster
\item"(i)"  $A\subset\Gamma$, and
\item"(ii)"  there is a section $L$ for $S$ 
in $\Lambda$ such that $A$ separates points on $L$ 
\RM(i.e., when $\ell,\ell'\in L$, $\ell\ne\ell'$, then
there is some $a\in A$ s.t. $e^{i2\pi\ell\cdot a}\ne 
e^{i2\pi\ell'\cdot a})$.
\endroster

Then it follows that every measurable section $D'$ inside 
$\Omega'$ 
\RM(reduced form\/\RM) for the action {\rm(3)} of $A$ by 
translation is a fundamental domain
for $\Gamma$ and, moreover, that
$$\Omega' =\bigcup_{a\in A/K} (D'+a)\tag 4$$
and
$$(D'+a_1)\cap (D'+a_2)=\varnothing$$
for all $a_1\ne a_2$ in $A/K$.
 \endproclaim

\subheading{\num{3.1.} Spectral duality}  
In studying more general simple factors, we introduced in 
[JP3] an {\it inductive limit construction} which
applies to the basic factors described in Theorem 1, and 
we found, as the limit object, the Hilbert space 
$L^2(\mu)$ where $\mu$ is a Hausdorff
measure of fractional dimension (see [Fa, Hu, St1--3]).  
Such measures are
known to be supported by Cantor type-sets, $\cC$, say
 (see [Hu]), but typically the
Lebesgue measure of $\cC$ is zero.  We now show that $\cC$ 
may be built by self-similarity from simple factors.

\subheading{\num{3.2}\rm}  Let $(\Omega ,\Lambda)$ be a 
spectral pair subject to the conditions in Theorem 1; let 
$K\subset A\subset\Gamma$ be the
associated lattices; let $L$ be the section in $\Lambda$ 
(assume $0\in L$); and
finally,
 let $R$ be the inclusion matrix for $K\subset\Gamma$.  (Let
$\{u_i\}^d_{i=1}$ be generators for $K$ over 
$\bz$ and $\{v_i\}^d_{i=1}$ for $\Gamma$; then $R\in 
M_d(\Bbb Z)\cap \roman{GL}_d(\Bbb R)$ may be defined by 
$u_i=\sum_j R_{ij}
v_j$.  Recall $K=\{\sum_i n_i u_i : n_i\in\bz\ ,\ 1\leq 
i\leq d\}$, and similarly for $\Gamma$.)  Since $L\subset 
K^0$, we may consider affine
mappings, $s\mapsto Rs+\ell$, acting on the lattice $K^0$. 
 This map will be denoted $\tau_\ell$, and the underlying 
lattice $K^0$ will be
understood from the context.  Consider the mapping 
$\tau_0(s)=Rs$ given by matrix-multiplication.  When the 
bases $(u_i)$
 for $K$ and $(v_i)$ for
$\Gamma$ are given, let $(u^*_i)$ for $K^0$ and $(v^*_i)$ 
for $\Gamma\,^0$ be dual bases, i.e., $u^*_i\cdot 
u_j=v^*_i\cdot v_j=\delta_{ij}$, $1\leq
i, j\leq d$.  For $s=\sum_i s_i u^*_i$ with integral 
coordinates, $s_i\in\bz$, we have
$$\tau_0(s)=\sum_i (Rs)_i u^*_i =\sum_i s_i v^*_i\ ,$$
and $(Rs)_i=\sum_j R_{ij} s_j$.  Note then that 
$\tau_0(K^0)\subset K^0$, and each $\tau_\ell$, $\ell\in 
L$, is affine on the lattice $K^0$.  If
$\Gamma\,^0$ is identified with a sublattice in $K^0$, 
then $\tau_0(K^0)=\Gamma\,^0$, and the matrix-transpose 
$R'_{ij}=R_{ji}$ {\it is} the
inclusion-matrix for the dual lattice-inclusion 
$\Gamma\,^0\subset K^0$.

\subheading{\num{3.3.} The Fractal Measure}  Also consider 
the affine maps $S_b$ on $\br^d$ given by
$$S_b x=R^{-1}x+b\ ,\qquad x\in\br^d\ .\tag 5$$
In formula (5), the term $R^{-1}x$ is really 
$\tau_0^{-1}(x)$, which is to say that the matrix-product 
$R^{-1}x$ must refer to the same 
basis $(u^*_i)$
(for $K^0$) that was used in calculating $\tau_0$ above.  
(In some different basis, of course, the matrix will 
change, i.e., $R$ becomes $ARA^{-1}$
with $A$ denoting the associated transform matrix.)

Let $N$ be the cardinality of $L$; by Theorem 1, it is 
also the order of the group $A/K$.  Pick a subset 
$\cB\subset A$, $0\in\cB$, representing
the elements in $A/K$, equivalently a section for the 
quotient; and let the affine maps $S_b$ be indexed by 
$b\in\cB$.  By Hutchinson's theorem (see
[Hu, St1--2]) there is self-similar probability measure 
$\mu$ on $\br^d$ such that $\mu={1\over N}\sum_{b\in\cB} 
\mu\circ S^{-1}_b\ ,$ or,
equivalently,
$$\int f(x)\,d\mu(x)={1\over N} \sum_{b\in\cB}\int 
f(S_bx)\, d\mu(x)$$
for measurable functions $f$ on $\br^d$.  We show in [JP3] 
that 
there is a ``Cantor set'' $\cC\subset\br^d$, which is 
built from iteration of the
decomposition (4) and self-similarity and which supports 
$\mu$, i.e., $\mu(\cC)=1$.  We 
let $L^2(\mu)$ be the corresponding  Hilbert space.

\subheading{\num{3.4.} The Cuntz Algebra}  
Our two theorems below connect the classical harmonic 
analysis of $(\Omega ,\Lambda)$ to the associated
fractal measure $\mu$:

\proclaim{Theorem 2}  Let $(\Omega ,\Lambda)$ be a 
nondegenerate simple factor given by the conditions in 
Theorem 
\RM1 with matrix $R$ for the lattice
inclusion $K\subset\Gamma$, and section $L$ for $\Lambda$ 
such that $0\in L$ and $\Lambda=L\dotplus\Gamma\,^0$, and 
finally let $\mu$ be the associated
Hutchinson measure with support $\cC$.  Then it follows that
\par
{\rm(i)}  
$\{ e_s:s\in K^0\}$ separates points in $\cC$, i.e., for 
$x\ne x'$ in $\cC$,
$\exists s\in K^0$ s.t. $e_s(x)\ne e_s(x')$.
\par
{\rm(ii)}  For each $\ell\in L$, an isometry $T_\ell$ 
acting on $L^2(\mu)$ is well defined by
$$T_\ell e_s=e_{\tau_{\ell}(s)}\quad \forall s\in K^0\ 
.\tag 6$$
\par
{\rm(iii)}  As operators on $L^2(\mu)$, the isometries 
$T_\ell$ satisfy
$$T^*_\ell T_{\ell'}=\cases 0&\text{if $\ell\ne\ell'$ in 
$L$},\\
I&\text{if $\ell=\ell'$,}\endcases\qquad\text{and}\qquad
\sum_{\ell\in L} T_\ell T^*_\ell = I$$
where $I$ denotes the identity operator on $L^2(\mu)$.
\par
{\rm(iv)}  The representation of the Cuntz $C^*$-algebra 
$\cO(L)$ generated by the isometries in 
{\rm(iii)} \RM(see {\rm[Cu])} has a canonical factor
decomposition associated with the triple $K\subset 
A\subset \Gamma$ of
lattices and the \RM(dual\/\RM) fractal measure $\mu$ may 
be reconstructed directly
from the associated factor state on $\cO(L)$ of the 
decomposition.  
\RM(Note that the decomposition is orthogonal, and in the 
category of
representations of C$^*$-algebras\/\RM; see {\rm[BR])}.
\par
{\rm(v)} The cyclic $e_o$-representation of $\cO(L)$ by 
the $T_\ell$ isometries is the GNS representation 
\RM(see {\rm[BR])} of the factor state $\omega$
on $\cO(L)$ which is determined by the relations in 
{\rm(iii)}, 
$\omega(I)=1$, and $\omega (T_oT^*_o)=1=\omega(T_0)$.
\par
{\rm(vi)} The set of all vectors
$$\{e_o\}\cup\bigcup_{n=1}^\infty \left\{ T_{\ell_1}\cdots 
T_{\ell_n} e_o : \ell_i\in L\right\}\tag 7$$
is maximal $\mu$-orthogonal and spans a closed subspace in 
$L^2(\mu)$ with infinite-dimen\-sion\-al orthogonal 
complement.
\par
{\rm(vii)} The Fourier transform
$$\widehat{\mu}(t) =\int_\cC e_t(x)\ d\mu(x)$$
satisfies the functional transformation law
$$\widehat{\mu}(Rt)=B(Rt)\widehat{\mu}(t)\quad \forall 
t\in\br^d\ ,$$
where
$$B(t)={1\over N} \sum_{b\in\cB} e^{i2\pi b\cdot t}$$
and $\widehat{\mu}(\cdot)$ has an associated infinite 
product-formula.
\endproclaim

\rem{Remark}  We view the representation (6) as a 
substitute for an orthogonal harmonic analysis for 
$L^2(\mu)$, with $\mu$ fractal, and note that
the relations in (iii) above have the flavor of an 
orthogonal double-decomposition but {\it not} an 
orthogonal expansion in the classical sense of
Fourier integrals (or series).  Indeed, Strichartz [St2] 
showed that there is not a direct way of making an exact 
classical Fourier decomposition for
$L^2(\mu)$ when $\mu$ is fractal.
\endrem

\heading 4.  Orthogonal frequencies in $L^2(\mu)$\endheading

Note that in (vi) the vectors from (7) are represented by 
orthogonal frequencies $e_\xi$ of the form (1) where $\xi$ 
is in the subset
$\cL(L)\subset\br^d$ of all affine sums (with $n$ variable):
$$\sum^n_{k=1} R^{k-1} 
\ell_k=\tau_{\ell_1}\tau_{\ell_2}\cdots\tau_{\ell_n}(0)\ 
,$$
$\forall\ell_k\in L$, and the $n=0$ term corresponding (by 
definition) to $\xi=0$.

\proclaim{Theorem 3 \rm(details [JP3])} {\rm(i)} $\{e_\xi 
:\xi\in\cL(L)\}$ is maximally orthogonal in $L^2(\mu)$.
\par{\rm(ii)}  None of the functions $e_s(x)=e^{i2\pi 
s\cdot x} (x\in\br^d)$ for $s\in\br^d\backslash\cL(L)$ is 
in the $L^2(\mu)$-closed linear span of
the pure frequencies of $\cL(L)$.  That is,
$$\sigma_L(s):=\sum_{\xi\in\cL(L)} 
|\widehat{\mu}(s-\xi)|^2 < 1$$
when $s$ is in $\br^d\backslash\cL(L)$.
\endproclaim

However, computer-calculations (Mathematica) show that
$$\sigma_L(s)=\Vert P_{\cL(L)} e_s\Vert^2_{L^2(\mu)}$$
is close to 1 (within third decimal place) 
when $s=(s_1,\dots ,s_d)\in K^0\backslash\cL(L)$ and 
$s_i>0$, $1\leq i\leq d$.

\heading 5.  Returning to $(\Omega ,\Lambda)$\endheading

Our final result shows that the system $(\Omega ,\Lambda)$ 
may be reconstructed from a given Cuntz-representation 
acting on $L^2(\mu)$.

\proclaim{Theorem 4}  Let $\mu$ be a probability measure 
on $\br^d$ with compact support, and let $K\subset\Gamma$ 
be a rank $d$ lattice system, with
inclusion matrix $R$.  Suppose a subset $L$ s.t. $0\in 
L\subset K^0$ induces operators $\{T_\ell\}_{\ell\in L}$ 
by 
{\rm(6)}, acting isometrically on
$L^2(\mu)$ and satisfying the Cuntz-relations {\rm(iii)} 
in Theorem \RM2.  
Then it follows that $\mu$ is a fractal measure which is 
generated by
self-similarity from some spectral pair $(\Omega 
,\Lambda)$ in $\br^d$ satisfying the conditions in Theorem 

\RM1 for nondegenerate simple factors.
\endproclaim

\heading Acknowledgment \endheading 
The authors gratefully acknowledge helpful correspondence 
from Professor R. S.
Strichartz on the connection of our earlier paper [JP1] to 
the harmonic
analysis of fractal measures. 


\Refs
\widestnumber\key{JP3}

\ref\key BR
\by O. Bratteli and D. W. Robinson
\book Operator algebras and quantum statistical mechanics, 
{\rm vol. I, revised
ed.}
\publ Springer-Verlag
\publaddr New York
\yr 1987
\endref

\ref\key Cu
\by J. Cuntz
\paper C$^*$-Algebras generated by isometries
\jour Comm. Math. Phys.
\vol 57
\yr 1977
\pages 173--185
\endref

\ref\key Fa
\by K. J. Falconer
\book The geometry of fractal sets
\publ Cambridge Univ. Press
\publaddr London and New York
\yr 1985
\endref

\ref\key Fu
\by B. Fuglede
\paper Commuting self-adjoint partial differential 
operators and a group theoretic problem
\jour J. Funct. Anal.
\vol 16
\yr 1974
\pages 101--121
\endref

\ref\key Hu
\by J. E. Hutchinson
\paper Fractals and self-similarity
\jour Indiana Univ. Math. J.
\vol 30
\yr 1981
\pages 713--747
\endref

\ref\key JP1
\by P. E. T. Jorgensen and S. Pedersen
\paper Spectral theory for Borel sets in $\br^n$ of finite 
measure
\jour J. Funct. Anal.
\vol 107
\yr 1992
\pages 72--104
\endref

\ref\key JP2
\bysame
\paper Sur un probl\`eme spectral alg\'ebrique
\jour C. R. Acad. Sci. Paris S\'er. I Math
\vol312
\yr 1991
\pages 495--498
\endref

\ref\key JP3
\bysame
\paper Spectral duality for C$^*$-algebras and fractal 
measures
\jour in preparation
\endref

\ref\key St1
\by R. S. Strichartz
\paper Self-similar measures and their Fourier transforms. 
{\rm I, II}
\jour Indiana Univ. Math. J.
\vol 39
\yr 1990
\pages 797--817; II, Trans. Amer. Math. Soc. {\bf336}
(1993), 335--361
\endref

\ref\key St2
\bysame
\paper Fourier asymptotics of fractal measures
\jour J. Funct. Anal.
\vol 89
\yr 1990
\pages 154--187
\endref

\ref\key St3
\bysame
\paper Wavelets and self-affine tilings
\jour Constr. Approx. 
\vol9 \yr 1993 \pages 327--346
\endref

\endRefs
\enddocument